\documentclass[12 pt, onecolumn]{IEEEtran}
\usepackage{amsmath,amssymb,euscript,latexsym,graphicx}
\usepackage{bbm,color,amstext,wasysym,subfig,parskip,balance}
\graphicspath{{./},{./figures/}}

\newtheorem{thm}{Theorem}

\newcommand{\be}{\begin{equation}}
\newcommand{\ee}{\end{equation}}
\newcommand{\mR}{{\mathbb R}}

\newcommand{\cH}{{\mathcal H}}

\newcommand{\cD}{{\mathcal D}}
\newcommand{\cE}{{\mathcal E}}

\newcommand{\cS}{{\mathcal S}}

\newcommand{\cU}{{\mathcal U}}
\newcommand{\cL}{{\mathcal L}}

\newcommand{\trace}{\operatorname{tr}}

\definecolor{grey}{rgb}{0.6,0.6,0.6}
\definecolor{lightgray}{rgb}{0.97,.99,0.99}

\begin{document}
\title{Interpolation of Matrices and Matrix-Valued Measures: The Unbalanced Case}
\author{Yongxin Chen, Tryphon T. Georgiou, and Allen Tannenbaum
\thanks{Y.\ Chen is with the Department of Medical Physics, Memorial Sloan Kettering Cancer Center, NY; email: chen2468@umn.edu}
\thanks{T.\ T. Georgiou is with the Department of Mechanical and Aerospace Engineering, University of California, Irvine, CA; email: tryphon@uci.edu}
\thanks{A.\ Tannenbaum is with the Departments of Computer Science and Applied Mathematics \& Statistics, Stony Brook University, NY; email: allen.tannenbaum@stonybrook.edu}}

\maketitle

\begin{abstract}
In this note, we propose an unbalanced version of the quantum mechanical version of optimal mass transport that was based on the Lindblad equation. We formulate a natural interpolation framework between density matrices and matrix-valued measures via a quantum mechanical formulation of Fisher-Rao information and the matricial Wasserstein distance.
\end{abstract}

\section{Introduction}
Optimal mass transport (OMT) besides its intrinsic mathematical elegance has proven to be a very powerful methodology for numerous problems in econometrics, systems and
control, information theory, statistical filtering and estimation, computer
vision, and signal/image processing \cite{Rachev,TGT}.
However, the standard requirement of optimal mass transport (OMT) of mass preservation is many times unnatural for real-world problems. For example, in image registration or optical flow, one must impose {\em ad hoc} normalizations on the imagery that do not have a physical justification. Accordingly, there have been a number of approaches based on the interpolation of the Wasserstein metric from OMT with some other metric such as $L^2$ \cite{Benamou,French} or some information-theoretic distance \cite{fisher1,hellinger}.
The seminal work of Benamou and Brenier \cite{French} makes this possible via certain modifications of the underlying energy functional and corresponding continuity equation.

In previous work, Chen {\em at al.} \cite{Chen} developed a quantum mechanical framework for defining a Wasserstein distance on matrix-valued densities (normalized to have trace 1), via a variational framework with a continuity equation constraint generalizing the  work of \cite{French}. We show in the present note, that the methodology and definitions
in \cite{Chen} of matrix analogues of the gradient and divergence, allow us to formulate in a rather straightforward manner a natural energy functional and continuity equation that generalize the information theoretic unbalanced approaches that give interpolations of Wasserstein and Fisher-Rao \cite{fisher1} and Hellinger \cite{hellinger}. Thus one derives unbalanced version of optimal mass transport in the matrix-valued case. We also show that the unbalanced problem may be formulated as one of convex optimization that makes it applicable to a variety of applications.

\section{Background on unbalanced mass transport}

In this section, we briefly introduce the basis of OMT and review two possible methods for an unbalanced version of OMT in the scalar case following \cite{Benamou,fisher1}. The original formulation of OMT is
	\begin{equation}
		\inf_{T} \left\{\int_{\mR^m} c(x, T(x))\rho_0(x)dx~\mid~ T_\sharp \rho_0 = \rho_1\right\},
	\end{equation}
where $c(x,y)$ denotes the cost of moving unit mass from $x$ to $y$, and $\rho_0, \rho_1$ are two given probability densities in $\mR^m$. The operator $(\cdot)_\sharp$ represents push forward. By relaxing the map $T$ to a coupling $\pi$ \cite{Kantorovich1948}, we obtain the Kantorovich form
	\begin{equation}
		\inf_{\pi \in \Pi(\rho_0,\rho_1)} \int_{\mR^m} c(x,y) \pi(dx,dy),
	\end{equation}
with $\Pi(\rho_0,\rho_1)$ denoting the set of all joint distributions (couplings) between $\rho_0$ and $\rho_1$. When the cost function $c(x,y)=\|x-y\|^2$, the problem has extremely rich structures. As firstly pointed out in \cite{French}, the OMT problem has the fluid dynamic formulation
	\begin{subequations}
	\begin{eqnarray}
	W_2(\rho_0,\rho_1)^2:=&&\inf_{\rho,v}  \int_0^1\int_{\mR^m} \rho(t,x) \|v(t,x)\|^2 \, dxdt
	\\&&\frac{\partial \rho}{\partial t} + \nabla \cdot (\rho v) = 0,
	\\&&\rho(0, \cdot) = \rho_0 (\cdot), \; \rho(1, \cdot) = \rho_1 (\cdot).
	\end{eqnarray}
	\end{subequations}
The optimal value gives the (squared) Wasserstein distance $W_2$. We next sketch two possible ways of extending $W_2$ to unbalanced measures. These formulation are used to study unbalance transport problem for matrices and matrix-valued densities.

\subsection{$L^2$ and OMT}

As noted in \cite{French,Benamou}, the $L^2$ problem can be used in conjunction with OMT in case of unbalanced mass distributions. The dynamic framework of \cite{French} makes this quite straightforward. Full details and numerics may be found in \cite{Benamou}.

Accordingly, given two unbalanced densities $\rho_0$ and $\rho_1$ it is natural to seek a distribution
$\tilde{\rho}_1$ the closest density to $\rho_1$ in the $L^2$ sense,
which minimizes the Wasserstein distance $W_2(\rho_0,
\tilde{\rho}_1)^2.$ The $L^2$ perturbation may be interpreted as ``noise.''
One can then show that this problem amounts to minimizing
\begin{subequations}
\begin{eqnarray}
&&\inf_{\rho,v,\tilde{\rho}_1}  \int_0^1\int_{\mR^m} \rho(t,x) \|v\|^2 \, dx \,
dt + \alpha \int_{\mR^m} (\rho_1(x) - \tilde{\rho}_1(x))^2 \,dx,
\\&&\frac{\partial \rho}{\partial t} + \nabla \cdot (\rho v) = 0,
\\&&\rho(0, \cdot) = \rho_0 (\cdot), \; \rho(1, \cdot) = \tilde{\rho}_1 (\cdot).
\end{eqnarray}
\end{subequations}
with  $\alpha >0$ being the coefficient balancing the two parts of the cost.
This method has been used in several applications including optical flow; see \cite{Mueller} and the references therein.

Here we introduce a slightly different form of interpolation distance between $W_2$ and $L^2$ that allows Riemannian structure. We bring in a source term $s$ in the continuity equation and construct a convex optimization problem
\begin{subequations}\label{eq:scalarunbalanceL}
\begin{eqnarray}
&&\inf_{\rho,v,s}  \int_0^1\int_{\mR^m}\left\{ \rho(t,x) \|v\|^2 +\alpha s(t,x)^2\right\}dx \,
dt,
\\&&\frac{\partial \rho}{\partial t} + \nabla \cdot (\rho v) = s, \label{eq:scalarunbalanceL2}
\\&&\rho(0, \cdot) = \rho_0 (\cdot), \; \rho(1, \cdot) = \rho_1 (\cdot).
\end{eqnarray}
\end{subequations}	
It can be shown that the square root of the minimum of the above is a well-defined metric on the space of probability densities with finite second order moments.

\subsection{Unbalanced mass transport: information theoretic formulations} \label{sec:unbalanced_info}

We now review how Fisher-Rao may be employed to get an unbalanced formulation of OMT \cite{fisher1}. Unlike the method of \cite{French, Benamou} where one interpolates $L^2$ and OMT, here the idea is to interpolate OMT and the Fisher-Rao metric. This is quite powerful since one explicitly combines an information theoretic method with Wasserstein, and thus makes contact with the recent work on Schroedinger bridges \cite{CGP}. There is a related method in \cite{hellinger} that interpolates between OMT and the Hellinger metric.

The interpolation of Fisher-Rao and Wasserstein is given as
\begin{subequations}\label{eq:scalarunbalancerao}
\begin{eqnarray}\label{eq:scalarunbalancerao1}
&&\inf_{\rho,v,r} \int_0^1\int _{\mR^m}\left\{ \rho(t,x) \|v\|^2 + \alpha \rho(t,x) r^2 \right\}dx \,dt
\\&&\frac{\partial \rho}{\partial t} + \nabla \cdot (\rho v) = \rho r \label{eq:scalarunbalancerao2}
\\&&
\rho(0, \cdot) = \rho_0 (\cdot), \; \rho(1, \cdot) = \rho_1 (\cdot) .
\end{eqnarray}
\end{subequations}
Here the minimum is taken over all the density flow $\rho$, velocity field $v$ and relative source intensity $r$ satisfying the continuity equation with source term \eqref{eq:scalarunbalancerao2}. Observing the source terms in \eqref{eq:scalarunbalanceL2} and \eqref{eq:scalarunbalancerao2} have the relation $s=\rho r$, we rewrite the second of the cost in \eqref{eq:scalarunbalancerao1} as
	\begin{equation}\label{eq:fisherrao}
		\int_0^1\int _{\mR^m}\rho(t,x) r(t,x)^2dx dt=\int_0^1\int _{\mR^m}\frac{s(t,x)^2}{\rho(t,x)}dx dt.
	\end{equation}
This should be compared to
	\begin{equation}\label{eq:L2}
	\int_0^1\int _{\mR^m}s(t,x)^2dxdt,
	\end{equation}
which is used in \eqref{eq:scalarunbalanceL}. The cost \eqref{eq:L2} corresponds to the $L^2$ metric while \eqref{eq:fisherrao} defines the Fisher-Rao between two smooth densities as
	\begin{eqnarray*}
	d_{FR}(\rho_0,\rho_1)^2:=&&\inf_{\rho,s}\int_0^1\int _{\mR^m}\frac{s(t,x)^2}{\rho(t,x)}dx dt
	\\&&\frac{\partial \rho}{\partial t}= s,
	\\&&\rho(0, \cdot) = \rho_0 (\cdot), \; \rho(1, \cdot) = \rho_1 (\cdot).
	\end{eqnarray*}

\section{Quantum continuity equation} \label{sec:background}

We sketch here the necessary background from \cite{Chen}. Consider two positive definite (Hermitian) matrices $\rho_0$ and $\rho_1$. We seek a suitable generalization of the continuity equation that links the two matrices with a smooth path within the cone of positive matrices in suitable ways. In the context of quantum mechanics, $\rho$ may represent a density matrix.
A standing assumption is that  $\trace(\rho_0)=\trace(\rho_1)=1$, and thereby, we seek paths $\rho(t)$ ($t\in[0,1]$) between the two that maintain the same value for the trace.

Let $\cH$ and $\cS$ denote the set of  $n\times n$ Hermitian and skew-Hermitian matrices, respectively. Since matrices are $n\times n$ throughout, we dispense of $n$ in the notation. We also denote the space of block-column vectors consisting of $N$ elements in $\cS$ and $\cH$ as $\cS^N$, respectively $\cH^N$. Let now let $\cH_+$ and $\cH_{++}$ denote the cones of nonnegative and positive definite matrices, respectively, and
	\[
	\cD_+ :=\{\rho \in \cH_{++} \mid \trace(\rho)=1\}.
	\]
Clearly, the tangent space of $\cD_+$, at any $\rho\in \cD_+$, is
	\[
	T_\rho=T:=\{ \sigma \in \cH \mid \trace(\sigma)=0\}.
	\]
We also use the standard notion of inner product
	\[
		\langle X, Y\rangle=\trace(X^*Y)
	\]
for both $\cH$ and $\cS$. For $X, Y\in \cH^N$ ($\cS^N$),
	\[
		\langle X, Y\rangle=\sum_{k=1}^N \trace(X_k^*Y_k).
	\]
Given $X=[X_1^*,\cdots,X_N^*]^* \in \cH^N$ ($\cS^N$), $Y\in \cH$ ($\cS$), denote
	\[
		XY=\left[\begin{array}{c}
		X_1\\
		\vdots \\
		X_N
		\end{array}\right]Y
		:=
		\left[\begin{array}{c}
		X_1Y\\
		\vdots \\
		X_N Y
		\end{array}\right],
	\]
and
	\[
		YX=Y\left[\begin{array}{c}
		X_1\\
		\vdots \\
		X_N
		\end{array}\right]
		:=
		\left[\begin{array}{c}
		YX_1\\
		\vdots \\
		YX_N
		\end{array}\right].
	\]

In a quantum system (open quantum system), the dynamics of density matrices can be described by the Lindblad equation
	\begin{equation}\label{eq:lindblad}
	\dot\rho =-i [H, \rho]+\sum_{k=1}^N( L_k\rho L_k^*-\frac{1}{2}\rho L_k^*L_k-\frac{1}{2}L_k^*L_k\rho).
	\end{equation}
Here the first term on the right hand side describes the evolution of the state under the effect of the Hamiltonian $H$ and it is energy preserving. The rest of the terms on the RHS represent the diffusion and capture the dissipation of energy. Note that this is the quantum analogue of the Laplacian operator $\Delta$.

In the following, assume $L_k=L_k^*$, i.e., $L_k\in \cH$ for all $k\in 1 \ldots, N$. Under this assumption, we can define
	\begin{equation}\label{eq:gradient}
		\nabla_L: \cH \rightarrow \cS^N, ~~X \mapsto
		\left[ \begin{array}{c}
		L_1 X-XL_1\\
		\vdots \\
		L_N X-X L_N
		\end{array}\right]
	\end{equation}
as the gradient operator. The dual of $\nabla_L$, which is an analogue of the divergence operator, is given by
	\begin{equation}\label{eq:divergence}
		\nabla_L^*: \cS^N \rightarrow \cH,~~Y=
		\left[ \begin{array}{c}
		Y_1\\
		\vdots \\
		Y_N
		\end{array}\right]
		\mapsto
		\sum_k^N L_k Y_k-Y_k L_k.
	\end{equation}
One can get this by definition
	\[
		\langle \nabla_L X, Y\rangle =\langle X, \nabla_L^* Y\rangle.
	\]
With this definition we calculate the ``Laplacian'' as
	\[
		\Delta_L X=-\nabla_L^*\nabla_L X=\sum_{k=1}^N( 2L_k\rho L_k^*-\rho L_k^*L_k-L_k^*L_k\rho),
	\]
which is exactly (with some scaling) the diffusion term in the Lindblad equation \eqref{eq:lindblad}. Therefore the Lindblad equation (under the assumption $L_k=L_k^*$) can be rewritten as
	\[
		\dot\rho =-i [H, \rho]+\frac{1}{2} \Delta_L \rho.
	\]
	
Note that the gradient operator $\nabla_L$ acts just like the standard gradient operator. Note that, in particular,
	\[
		\nabla_L(XY+YX)=\nabla_L X Y+X\nabla_L Y+\nabla_L YX+Y\nabla_L X,~~\forall X, Y\in \cH.
	\]
Using this gradient operator \eqref{eq:gradient}, we can then come up with several notions of the continuity equation. In the present note for the interpolation of Fisher-Rao and Wasserstein, we will use
	\begin{equation}\label{eq:continuity}
		\dot \rho=\frac{1}{2} \nabla_L^* (\rho v+v\rho),
	\end{equation}
where the ``velocity'' field $v=[v_1^*,\ldots,v_N^*]^*\in \cS^N$. Note $\rho v+v\rho \in \cS^N$, which is consistent with the definition of $\nabla_L^*$.

Usually, in the Lindblad equation \eqref{eq:lindblad}, $N$ is taken to be $n^2-1$. However, in general, we may choose $N\le n^2-1$, as needed, possibly large enough such that in \eqref{eq:continuity} we are able to cover the whole tangent space $T_\rho$ at $\rho$ for all $\rho\in\cD_+$.
In particular, we need $\nabla_L$ to have the property that the identity matrix $I$ spans its null space. For instance, one can choose $L_1,\ldots,L_N$ to be a basis of the Hermitian matrices $\cH$, in which case $N=n(n+1)/2$. Obviously this construction ensures that the null space of $\nabla_L$ is spanned by $I$.

 Now we have everything ready to define the fluid dynamic formulation of optimal transport on the space $\cD_+$ of density matrices. Given two density matrices $\rho_0, \rho_1 \in \cD_+$, one can formulate the following optimization problem
 	\begin{subequations}\label{eq:quantumomt}
 	\begin{eqnarray}\label{eq:quantumomt1}
	W_2(\rho_0, \rho_1)^2:=&&\inf_{\rho\in \cD_+, v\in \cS^N} \int_0^1 \trace(\rho v^*v) dt\\
	&&\dot \rho=\frac{1}{2} \nabla_L^* (\rho v+v\rho), \label{eq:quantumomt2}\\
	&& \rho(0)=\rho_0, ~~\rho(1)=\rho_1  \label{eq:quantumomt3}
	\end{eqnarray}
	\end{subequations}
and define the ``Wasserstein distance'' between $\rho_0$ and $\rho_1$ to be the square root of the minimum of the cost \eqref{eq:quantumomt1}. Note here for $v\in \cS^N$, $v^*v=\sum_{k=1}^N v_k^*v_k$. The Wasserstein distance function $ W_2(\rho, \rho+\delta\rho)$ gives an Riemannian structure on the tangent space $T_\rho$, and therefore $W_2(\cdot, \cdot)$ indeed defines a metric on $\cD_+$. One can introduce a Lagrangian multiplier $\lambda \in \cH$ for the constraints \eqref{eq:quantumomt2} and arrive at the following sufficient conditions for optimality.
\begin{thm}
Suppose there exists $\lambda(\cdot)\in \cH$ satisfying
	\begin{subequations}
	\begin{equation}
		\dot\lambda=\frac{1}{2}(\nabla_L \lambda)^* (\nabla_L \lambda)
		=\frac{1}{2}\sum_{k=1}^N(\nabla_L \lambda)_k^* (\nabla_L \lambda)_k
	\end{equation}
	such that the solution of
	\begin{equation}
		\dot \rho=-\frac{1}{2} \nabla_L^* (\rho \nabla_L \lambda+\nabla_L \lambda\rho)
	\end{equation}
	\end{subequations}
	matches the two marginals $\rho(0)=\rho_0, \rho(1)=\rho_1$, then $(\rho, v=-\nabla_L \lambda)$ solves \eqref{eq:quantumomt}.
\end{thm}

\section{Interpolation of matrices: unbalanced case} \label{sec:unbalanced}
In this section, we formulate the main result of the present note, namely the  interpolation between quantum Wasserstein and Fisher-Rao and that between quantum Wasserstein and Frobenius norm, as generalizations of \eqref{eq:scalarunbalancerao} and \eqref{eq:scalarunbalanceL}, repectively.

\subsection{Interpolation between Wasserstein and Fisher-Rao}
Given $\rho_0,\rho_1\in \cH_{++}$ and $\alpha>0$, define
 	\begin{subequations}\label{eq:quantumomtunb}
 	\begin{eqnarray}\label{eq:quantumomtunb1}
	W_{2,FS}(\rho_0,\rho_1)^2:=&&\inf_{\rho\in \cH_{++}, v\in \cS^N, r\in\cH} \int_0^1\{ \trace(\rho v^*v)+
	\alpha\trace(\rho r^2)\} dt\\
	&&\dot \rho=\frac{1}{2} \nabla_L^* (\rho v+v\rho)+\frac12(\rho r+r \rho), \label{eq:quantumomtunb2}\\
	&& \rho(0)=\rho_0, ~~\rho(1)=\rho_1.  \label{eq:quantumomtunb3}
	\end{eqnarray}
	\end{subequations}
Note here the ``continuity'' equation \eqref{eq:quantumomtunb2}, as a non-commutative generalization of \eqref{eq:scalarunbalancerao}, preserves positivity but \textbf{\emph{not}} mass. This distance $V_2$ is an interpolation of $W_2$ and the Fisher-Rao distance
	\begin{eqnarray*}
	d_{FR}(\rho_0,\rho_1):= &&\inf_{r\in\cH} \trace(\rho r^2) dt
	\\&& \dot \rho=\frac12(\rho r+r \rho)
	\\&&\rho(0)=\rho_0, ~~\rho(1)=\rho_1.
	\end{eqnarray*}
Recall that the Bures metric \cite{Uhl92} on the space of density matrices $\cD$ is defined as
	\[
		d_B(\rho,\rho+\delta\rho)^2=\frac12\trace(G\delta\rho)
	\]
where $G\in\cH$ is the unique solution of
	\[
		\rho G+G\rho=\delta\rho.
	\]
on $\cH$. It follows
	\[
		\frac12\trace(G\delta\rho)=\trace(\rho G^2).
	\]
Hence, $d_{FS}$ is equivalent to the Bures metric $d_B$ when restricted to $\cD$.

We next exam the optimality condition for \eqref{eq:quantumomtunb}. Let $\lambda(\cdot) \in \cH$ be a smooth Lagrangian multiplier  for the constraints \eqref{eq:quantumomtunb2} and construct the Lagrangian
	\begin{eqnarray*}
		\cL(\rho,v,r,\lambda)&=&\int_0^1 \left\{\frac{1}{2}\trace(\rho v^*v)+\frac{\alpha}{2}\trace(\rho r^*r)-\trace(\lambda
		(\dot \rho-\frac{1}{2} \nabla_L^* (\rho v+v\rho)-\frac12(\rho r+r \rho)))\right\}dt
		\\&=& \int_0^1 \left\{\frac{1}{2}\trace(\rho v^*v)+\frac{1}{2}\trace((\nabla_L\lambda)^* (\rho v+v\rho))
		+\frac{\alpha}{2}\trace(\rho r^*r)+\frac12\trace(\lambda(\rho r+r \rho))
		+\trace(\dot\lambda\rho)\right\}dt
		\\&&-\trace(\lambda(1)\rho_1)+\trace(\lambda(0)\rho_0).
	\end{eqnarray*}
Point-wise minimizing the above over $v$ yields
	\[
		v_{opt}(t)=-\nabla_L \lambda(t),
	\]
and
	\[
		r_{opt}(t)=-\frac{1}{\alpha}\lambda(t).
	\]
The corresponding minimum is
	\[
		\int_0^1 \left\{-\frac{1}{2}\trace(\rho (\nabla_L \lambda)^* (\nabla_L \lambda))-\frac{1}{2\alpha}\trace(\rho \lambda^2)+
		\trace(\dot\lambda\rho)\right\}dt -\trace(\lambda(1)\rho_1)+\trace(\lambda(0)\rho_0),
	\]	
from which we conclude the following sufficient conditions for optimality.
\begin{thm}\label{thm:matrixomtoptunb}
Suppose there exists $\lambda(\cdot)\in \cH$ satisfying
	\begin{subequations}
	\begin{equation}
		\dot\lambda=\frac{1}{2}(\nabla_L \lambda)^* (\nabla_L \lambda)+\frac{1}{2\alpha}\lambda^2
	\end{equation}
	such that the solution of
	\begin{equation}
		\dot \rho=-\frac{1}{2} \nabla_L^* (\rho \nabla_L \lambda+\nabla_L \lambda\rho)
		-\frac{1}{2\alpha}(\rho\lambda+\lambda\rho)
	\end{equation}
	\end{subequations}
	matches the marginals $\rho(0)=\rho_0, \rho(1)=\rho_1$. Then the triple $(\rho, v=-\nabla_L \lambda, r=-\frac{1}{\alpha}\lambda)$ solves \eqref{eq:quantumomtunb}.
\end{thm}

As in the classical OMT \cite{Villani}, the new distance $W_{2,FS}$ defines a ``Riemannian'' structure on  $\cH_{++}$. Given two tangent vectors $\delta_1, \delta_2\in \cH$ at $\rho$, the Riemannian metric is
	\begin{equation}\label{eq:riemannmetric}
		\langle \delta_1,\delta_2 \rangle_\rho = \frac12\trace(\rho\nabla\lambda_1^* \nabla\lambda_2
		+\rho\nabla\lambda_2^* \nabla\lambda_1)
		+\frac{1}{2\alpha}\trace(\rho\lambda_1\lambda_2+\rho\lambda_2\lambda_1),
	\end{equation}
where $\lambda_i\in \cH, ~(i=1,2)$ is the unique solution to
	\begin{equation}\label{eq:riemann}
		\delta_i=-\frac{1}{2} \nabla_L^* (\rho \nabla_L \lambda_i+\nabla_L \lambda_i\rho)
		-\frac{1}{2\alpha}(\rho\lambda_i+\lambda_i\rho).
	\end{equation}
In fact, $(v=-\nabla_L \lambda_i, r=-\frac{1}{\alpha}\lambda_i)$, with $\lambda_i$ being the solution  to the above equation \eqref{eq:riemann}, is the unique minimizer of
	\[
		\trace(\rho v^*v)+\alpha\trace(\rho r^2)
	\]
over all the $(v\in\cS^N, r\in \cH)$ satisfying
	\[
		\delta_i=\frac{1}{2} \nabla_L^* (\rho v+v\rho)+\frac12(\rho r+r \rho).
	\]	
It can be shown that our distance $W_{2,FS}$ is the geodesic distance on $\cH_{++}$ with Riemannian metric \eqref{eq:riemannmetric}, therefore, we have
	\[
		W_{2,FS}(\rho_0,\rho_1)=\inf_{\rho} \int_0^1 \sqrt{\langle \dot\rho(t),\dot\rho(t)\rangle_{\rho(t)}} dt.
	\]
In addition, the solution $\rho(\cdot)$ of \eqref{eq:quantumomtunb} possesses the nice property
	\[
		W_{2,FS}(\rho(s),\rho(t))=(t-s) W_{2,FS}(\rho_0,\rho_1)
	\]
for all $0\le s<t\le 1$. We remark that even though the Riemannian metric \eqref{eq:riemannmetric} is well defined on the boundary of $\cH_+$, $W_{2,FS}$ can be extended to $\cH_+$ by continuity with little effort.

As in the balanced case \cite{Chen}, \eqref{eq:quantumomtunb} has the following convex reformulation
	\begin{subequations}\label{eq:quantumomtunbcvx}
 	\begin{eqnarray}\label{eq:quantumomtunbcvx1}
	&&\inf_{\rho, u,s} \int_0^1 \{\trace(u^*\rho^{-1}u)+\alpha\trace(s^*\rho^{-1}s) \}dt,\\
	&&\dot \rho=\frac{1}{2} \nabla_L^* (u-\bar u)+\frac{1}{2}(s+s^*), \label{eq:quantumomtunbcvx2}\\
	&& \rho(0)=\rho_0, ~~\rho(1)=\rho_1  \label{eq:quantumomtunbcvx3}.
	\end{eqnarray}
	\end{subequations}
Here we simply used the change of variables $u=\rho v$ and $s=\rho r$. This convex formulation makes $W_{2,FS}$ suitable for various applications.

\subsection{Interpolation between Wasserstein and Frobenius}
As a straightforward generalization of \eqref{eq:scalarunbalanceL}, we define, for $\rho_0, \rho_1\in \cH_{++}$,
	\begin{subequations}\label{eq:quantumomtunbF}
	\begin{eqnarray}\label{eq:quantumomtunbF1}
	W_{2,F}(\rho_0,\rho_1)^2:=&&\inf_{\rho\in \cH_{++}, v\in \cS^N, s\in\cH} \int_0^1\{ \trace(\rho v^*v)+
	\alpha\trace(s^2)\} dt\\
	&&\dot \rho=\frac{1}{2} \nabla_L^* (\rho v+v\rho)+s, \label{eq:quantumomtunbF2}\\
	&& \rho(0)=\rho_0, ~~\rho(1)=\rho_1.  \label{eq:quantumomtunbF3}
	\end{eqnarray}
	\end{subequations}
The second part of the cost corresponds to the Frobenius metric. More specifically, the Frobenius metric can be rewritten as
	\[
		\|\rho_0-\rho_1\|_F^2=\inf_{s} \left\{\int_0^1 \trace(s^2) dt~\mid~\dot\rho=s, \rho(0)=\rho_0, \rho(1)=\rho_1\right\}.
	\]
Employing a similar Lagrangian argument, we obtain the optimality condition as follows.
\begin{thm}\label{thm:matrixomtoptunbF}
Suppose there exists $\lambda(\cdot)\in \cH$ satisfying
	\begin{subequations}
	\begin{equation}
		\dot\lambda=\frac{1}{2}(\nabla_L \lambda)^* (\nabla_L \lambda)
	\end{equation}
	such that the solution of
	\begin{equation}
		\dot \rho=-\frac{1}{2} \nabla_L^* (\rho \nabla_L \lambda+\nabla_L \lambda\rho)
		-\frac{1}{\alpha}\lambda
	\end{equation}
	\end{subequations}
	matches the marginals $\rho(0)=\rho_0, \rho(1)=\rho_1$. Then the triple $(\rho, v=-\nabla_L \lambda, s=-\frac{1}{\alpha}\lambda)$ solves \eqref{eq:quantumomtunbF}.
\end{thm}

Clearly, given any two tangent vectors $\delta_1, \delta_2$ at $\rho$ on $\cH_{++}$, the inner product
	\begin{equation}\label{eq:riemannmetricF}
	\langle \delta_1,\delta_2 \rangle_\rho = \frac12\trace(\rho\nabla\lambda_1^* \nabla\lambda_2
	+\rho\nabla\lambda_2^* \nabla\lambda_1)
		+\frac{1}{\alpha}\trace(\lambda_1\lambda_2),
	\end{equation}
endows the manifold $\cH_{++}$ an ``Riemannian'' like structure. Here $\lambda_i\in \cH, ~(i=1,2)$ is the unique solution to
	\begin{equation}\label{eq:riemannF}
		\delta_i=-\frac{1}{2} \nabla_L^* (\rho \nabla_L \lambda_i+\nabla_L \lambda_i\rho)
		-\frac{1}{\alpha}\lambda_i.
	\end{equation}
Besides, $(v=-\nabla_L \lambda_i, s=-\frac{1}{\alpha}{\lambda_i})$ minimizes
	\[
		\trace(\rho v^*v)+\alpha\trace(s^2)
	\]
over all the pairs $(v\in\cS^N, s\in \cH)$ satisfying
	\[
		\delta_i=\frac{1}{2} \nabla_L^* (\rho v+v\rho)+s.
	\]	
Similar to $W_{2,FS}$, it can be shown that our distance $W_{2,F}$ is the geodesic distance on $\cH_{++}$ with Riemannian metric \eqref{eq:riemannmetricF}, therefore, we have
	\[
		W_{2,F}(\rho_0,\rho_1)=\inf_{\rho} \int_0^1 \sqrt{\langle \dot\rho(t),\dot\rho(t)\rangle_{\rho(t)}} dt.
	\]
In addition, the solution $\rho(\cdot)$ of \eqref{eq:quantumomtunbF} possesses the nice property
	\[
		W_{2,F}(\rho(s),\rho(t))=(t-s) W_{2,F}(\rho_0,\rho_1)
	\]
for all $0\le s<t\le 1$.
Even though the Riemannian metric \eqref{eq:riemannmetricF} is well defined on the boundary of $\cH_+$, $W_{2,F}$ can be extended to $\cH_+$ by continuity.

Again, through changing of variable $u=\rho v$, we reformulate \eqref{eq:quantumomtunbF} as the following convex optimization problem
	\begin{subequations}\label{eq:quantumomtunbFcvx}
 	\begin{eqnarray}\label{eq:quantumomtunbFcvx1}
	&&\inf_{\rho, u,s} \int_0^1 \{\trace(u^*\rho^{-1}u)+\alpha\trace(s^2) \}dt,\\
	&&\dot \rho=\frac{1}{2} \nabla_L^* (u-\bar u)+s, \label{eq:quantumomtunbFcvx2}\\
	&& \rho(0)=\rho_0, ~~\rho(1)=\rho_1  \label{eq:quantumomtunbFcvx3}.
	\end{eqnarray}
	\end{subequations}

\section{Interpolation of matrix-valued measures: unbalanced case}

In applications it is often the case that one has to deal with matrix-valued distributions on dimensions which may represent space or frequency. Thus, in this case, the $\rho$'s may be $\cH_+$-valued functions on $E\subset\mR^m$. For instance, in the context of multivariable time series analysis it is natural to consider $m=1$; see, e.g., \cite{Lipeng}.
For simplicity, we assume $E$ to be a (convex) connected compact set.
Therefore, in this section
\begin{equation}
\cE=\{ \rho(\cdot) \mid \rho(x)\in \cH_+ \mbox{ for }x\in E \mbox{ such that }\int_{E} \trace(\rho(x))dx <\infty\}.
\end{equation}
Let $\cE_+$ denote the interior of $\cE$. Note that the problem on the subspace $\int_{E} \trace(\rho(x))dx=1$, i.e., the balanced case, has been studied in \cite{Chen}.
By combining the standard continuity equation on the Euclidean space and the continuity equation for positive definite matrices \eqref{eq:continuity}, and taking into the sources term, we obtain a continuity equation on $\cE_+$ for the flow $\rho(t,x)$ as
	\begin{equation}\label{eq:continuityspec}
		\frac{\partial \rho}{\partial t}+\frac{1}{2}\nabla_x\cdot(\rho w+w\rho)-\frac{1}{2} \nabla_L^* (\rho v+v\rho)
		-\frac12(\rho r+r \rho)=0,
	\end{equation}
or simply
	\begin{equation}\label{eq:continuityspecF}
		\frac{\partial \rho}{\partial t}+\frac{1}{2}\nabla_x\cdot(\rho w+w\rho)-\frac{1}{2} \nabla_L^* (\rho v+v\rho)
		-s=0.
	\end{equation}
Here $\nabla_x \cdot$ is the standard divergence operator on $\mR^m$, $w(t,x)\in {\cH}^m$ is the velocity field along the space dimension, and $v(t,x)\in \cS^N$ is the quantum velocity as before. We next present, based on the continuity equations, both the interpolating distance between Wasserstein and Fisher-Rao, and the interpolating distance between Wasserstein and Frobenius metric.

\subsection{Interpolation between Wasserstein and Fisher-Rao}
A dynamic formulation of matrix-valued optimal mass transport between
two given marginals $\rho_0, \rho_1 \in \cE_+$ ensues, namely,
	 \begin{subequations}\label{eq:matrixomt}
 	\begin{eqnarray}\label{eq:matrixomt1}
	W_{2,FS}(\rho_0, \rho_1)^2:=&&\inf_{\rho\in \cE_+, w\in {\cH}^m, v\in \cS^N,r\in\cH} \int_0^1\int_{\mR^m}
	\left \{\trace(\rho w^*w)+\gamma\trace(\rho v^*v)+\alpha\trace(\rho r^2)\right\}dxdt\\
	&&\frac{\partial \rho}{\partial t}+\frac{1}{2}\nabla_x\cdot(\rho w+w\rho)-\frac{1}{2} \nabla_L^* (\rho v+v\rho)
	-\frac12(\rho r+r \rho)=0, \label{eq:matrixomt2}\\
	&& \rho(0,\cdot)=\rho_0, ~~\rho(1,\cdot)=\rho_1  \label{eq:matrixomt3}.
	\end{eqnarray}
	\end{subequations}
The coefficient $\gamma>0$ is arbitrary and weighs in the relative significance of the two velocity fields. We then define the interpolating distance $W_{2,FS}(\rho_0, \rho_1)$ between $\rho_0$ and $\rho_1$ via \eqref{eq:matrixomt1}.

A sufficient condition for optimality can be obtained in a similar manner as before.
Here, we let $\lambda(\cdot,\cdot)\in \cH$ be a smooth function and define the Lagrangian
	\begin{align*}
		\cL(\rho,v,w, \lambda)&=
		\int_0^1\int_{\mR^m}\left \{  \frac{1}{2}\trace(\rho w^*w)+\frac{\gamma}{2}\trace(\rho v^*v)
		+\frac{\alpha}{2}\trace(\rho r^2)\right.\\
		&\left.
		-\trace(\lambda(\frac{\partial \rho}{\partial t}+\frac{1}{2}\nabla_x\cdot(\rho w+w\rho)-\frac{1}{2} \nabla_L^* (\rho v+v\rho)-\frac12(\rho r+r \rho)))\right\}dxdt.
	\end{align*}
Integration by parts yields
	\begin{eqnarray*}
	&&\int_0^1\int_{\mR^m}\left \{  \frac{1}{2}\trace(\rho w^*w)+\frac{\gamma}{2}\trace(\rho v^*v)
	+\frac{\alpha}{2}\trace(\rho r^2)\right.
	\\&&\left.+\trace(\frac{\partial \lambda}{\partial t}\rho)+\frac{1}{2}\langle\nabla_x\lambda,\rho w+w\rho\rangle+\frac{1}{2}
		\langle \nabla_L\lambda, \rho v+v\rho\rangle+\frac12\trace(\lambda(\rho r+r \rho))\right\}dxdt
	\end{eqnarray*}
Here we have discarded the terms on $\rho_0, \rho_1$. Minimizing the above pointwise over $w, v$ gives expressions for the optimal values as
	\[
		w_{opt}(t,x)=-\nabla_x \lambda(t,x),
	\]
	\[
		v_{opt}(t,x)=-\frac{1}{\gamma}\nabla_L \lambda(t,x),
	\]
and
	\[
		r_{opt}(t,x)=-\frac{1}{\alpha}\lambda(t,x).
	\]
Substituting these back to the Lagrangian we obtain
	\[
		\int_0^1\int_{\mR^m}\left \{  -\frac{1}{2}\trace(\rho (\nabla_x\lambda)^*(\nabla_x\lambda))
		-\frac{1}{2\gamma}\trace(\rho (\nabla_L \lambda)^*(\nabla_L\lambda))
		-\frac{1}{2\alpha}\trace(\rho \lambda^2)
		+\trace(\rho\frac{\partial \lambda}{\partial t})\right\}dxdt,
	\]	
and the sufficient conditions for optimality given below follow.
\begin{thm}
Suppose there exists smooth $\lambda(\cdot,\cdot)\in \cH$ satisfying
	\begin{subequations}
	\begin{equation}
		\frac{\partial\lambda}{\partial t}-\frac{1}{2}(\nabla_x\lambda)^*(\nabla_x\lambda)
		-\frac{1}{2\gamma}(\nabla_L \lambda)^* (\nabla_L \lambda)-\frac{1}{2\alpha}\lambda^2=0
	\end{equation}
such that the solution of
	\begin{equation}
		\frac{\partial \rho}{\partial t}-\frac{1}{2}\nabla_x\cdot(\rho \nabla_x \lambda+\nabla_x \lambda\rho)+
		\frac{1}{2\gamma} \nabla_L^* (\rho\nabla_L \lambda+\nabla_L \lambda\rho)
		+\frac{1}{2\alpha}(\rho\lambda+\lambda\rho)=0
	\end{equation}
	\end{subequations}
	matches the two marginals $\rho(0,\cdot)=\rho_0, \rho(1,\cdot)=\rho_1.$ Then $(\rho, w=-\nabla_x \lambda, v=-\frac{1}{\gamma}\nabla_L \lambda, r=-\frac{1}{\alpha}\lambda)$ solves \eqref{eq:matrixomt}.
\end{thm}
%
%
%

The distance $W_{2,FS}$ induces an Riemannian structure on $\cE_+$ and on the top of that, $W_{2,FS}$ is the corresponding geodesic distance. Since the discussion is similar to that in Section \ref{sec:unbalanced}, we skip the details here. As noted earlier, \eqref{eq:matrixomt} can again be cast as a convex optimization problem
 \begin{subequations}\label{eq:matrixomtcvx}
 	\begin{eqnarray}\label{eq:matrixomtcvx1}
	&&\inf_{\rho, q, u,s} \int_0^1\int_{\mR^m}
	\left \{\trace(q^*\rho^{-1}q)+\gamma\trace(u^*\rho^{-1}u)+\alpha\trace(s^*\rho^{-1}s)\right\}dxdt\\
	&&\frac{\partial \rho}{\partial t}+\frac{1}{2}\nabla_x\cdot(q+\bar{q})-\frac{1}{2} \nabla_L^* (u-\bar{u})-\frac{1}{2}(s+s^*)=0, \label{eq:matrixomtcvx2}\\
	&& \rho(0,\cdot)=\rho_0, ~~\rho(1,\cdot)=\rho_1  \label{eq:matrixomtcvx3}.
	\end{eqnarray}
	\end{subequations}

\subsection{Interpolation between Wasserstein and Frobenius}	
Given $\rho_0,\rho_1\in \cE_+$, we define the interpolating distance between Wasserstein and Frobenius through
	\begin{subequations}\label{eq:matrixomtF}
 	\begin{eqnarray}\label{eq:matrixomtF1}
	W_{2,F}(\rho_0, \rho_1)^2:=&&\inf_{\rho\in \cE_+, w\in {\cH}^m, v\in \cS^N,s\in\cH} \int_0^1\int_{\mR^m}
	\left \{\trace(\rho w^*w)+\gamma\trace(\rho v^*v)+\alpha\trace(s^2)\right\}dxdt\\
	&&\frac{\partial \rho}{\partial t}+\frac{1}{2}\nabla_x\cdot(\rho w+w\rho)-\frac{1}{2} \nabla_L^* (\rho v+v\rho)
	-s=0, \label{eq:matrixomtF2}\\
	&& \rho(0,\cdot)=\rho_0, ~~\rho(1,\cdot)=\rho_1  \label{eq:matrixomtF3}.
	\end{eqnarray}
	\end{subequations}
The optimality condition can be again established using Lagrangian method.
\begin{thm}
Suppose there exists smooth $\lambda(\cdot,\cdot)\in \cH$ satisfying
	\begin{subequations}
	\begin{equation}
		\frac{\partial\lambda}{\partial t}-\frac{1}{2}(\nabla_x\lambda)^*(\nabla_x\lambda)
		-\frac{1}{2\gamma}(\nabla_L \lambda)^* (\nabla_L \lambda)=0
	\end{equation}
such that the solution of
	\begin{equation}
		\frac{\partial \rho}{\partial t}-\frac{1}{2}\nabla_x\cdot(\rho \nabla_x \lambda+\nabla_x \lambda\rho)+
		\frac{1}{2\gamma} \nabla_L^* (\rho\nabla_L \lambda+\nabla_L \lambda\rho)
		+\frac{1}{\alpha}\lambda=0
	\end{equation}
	\end{subequations}
	matches the two marginals $\rho(0,\cdot)=\rho_0, \rho(1,\cdot)=\rho_1.$ Then $(\rho, w=-\nabla_x \lambda, v=-\frac{1}{\gamma}\nabla_L \lambda, s=-\frac{1}{\alpha}\lambda)$ solves \eqref{eq:matrixomtF}.
\end{thm}
Substituting $q=\rho w$ and $u=\rho v$, we get a convex formulation
\begin{subequations}\label{eq:matrixomtcvxF}
 	\begin{eqnarray}\label{eq:matrixomtcvxF1}
	&&\inf_{\rho, q, u,s} \int_0^1\int_{\mR^m}
	\left \{\trace(q^*\rho^{-1}q)+\gamma\trace(u^*\rho^{-1}u)+\alpha\trace(s^2)\right\}dxdt\\
	&&\frac{\partial \rho}{\partial t}+\frac{1}{2}\nabla_x\cdot(q+\bar{q})-\frac{1}{2} \nabla_L^* (u-\bar{u})-s=0, \label{eq:matrixomtcvxF2}\\
	&& \rho(0,\cdot)=\rho_0, ~~\rho(1,\cdot)=\rho_1  \label{eq:matrixomtcvxF3}.
	\end{eqnarray}
	\end{subequations}
	
\section{Gradient flow}
For completeness, we now derive some results for the gradient flow on $\cH_{++}$ of some energy functions with respect to $W_{2,FS}$ and $W_{2,F}$. We consider two functions: i) $\cS(\rho)=-\trace(\rho\log\rho-\rho)$; ii) $\cU(\rho)=\frac12\trace[(\rho-\hat\rho)^2]$.

\subsection{Gradient flow with respect to $W_{2,FS}$}
Taking the derivative of $\cS(\rho(t))$ over time yields
	\begin{eqnarray*}
	\frac{d\cS(\rho(t))}{dt} &=& -\trace(\log\rho \dot\rho)
	\\&=& -\trace(\log\rho(\frac{1}{2} \nabla_L^* (\rho v+v\rho)+\frac12(\rho r+r \rho)))
	\\&=& -\frac12\trace((\nabla_L\log\rho)^*(\rho v+v\rho)+\log\rho(\rho r+r\rho))
	\\&=& -\frac12\trace(\rho((\nabla_L\log\rho)^*v+v^*\nabla_L\log\rho)+\rho(\log\rho\cdot r+r\log\rho)).
	\end{eqnarray*}
This, together with the Riemannian metric \eqref{eq:riemannmetric}, points to the steepest ascent direction
	 \[
	 	v=-\nabla_L \log\rho,~~ r=-\frac{1}{\alpha} \log\rho.
	 \]
Hence, the gradient flow of $\cS$ with respect to $W_{2,FS}$ is
	\[
		\dot\rho(t)=-\frac{1}{2} \nabla_L^* (\rho \nabla_L\log\rho+\nabla_L\log\rho\cdot\rho)
		-\frac{1}{2\alpha}(\rho \log\rho+\log\rho\cdot \rho).
	\]	
Similarly, for $\cU$ we have
	\begin{eqnarray*}
	\frac{d\cU(\rho(t))}{dt}&=& \trace((\rho-\hat\rho)\dot\rho)
	\\&=&\trace((\rho-\hat\rho)(\frac{1}{2} \nabla_L^* (\rho v+v\rho)+\frac12(\rho r+r \rho)))
	\\&=&\frac12\trace(\rho((\nabla_L(\rho-\hat\rho))^*v+v^*\nabla_L(\rho-\hat\rho))+\rho((\rho-\hat\rho) r+r(\rho-\hat\rho))),
	\end{eqnarray*}
which leads to the steepest descent direction
	\[
		v=-\nabla_L(\rho-\hat\rho),~~r=-\frac{1}{\alpha}(\rho-\hat\rho),
	\]
and the gradient flow
	\[
		\dot\rho(t)=-\frac{1}{2} \nabla_L^* (\rho \nabla_L(\rho-\hat\rho)+\nabla_L(\rho-\hat\rho)\rho)
		-\frac{1}{2\alpha}(\rho (\rho-\hat\rho)+(\rho-\hat\rho) \rho).
	\]
	
\subsection{Gradient flow with respect to $W_{2,F}$}
The derivative of $\cS(\rho(t))$ over time is
	\begin{eqnarray*}
	\frac{d\cS(\rho(t))}{dt} &=& -\trace(\log\rho \dot\rho)
	\\&=& -\trace(\log\rho(\frac{1}{2} \nabla_L^* (\rho v+v\rho)+s))
	\\&=& -\trace(\frac12\rho((\nabla_L\log\rho)^*v+v^*\nabla_L\log\rho)+\log\rho\cdot s).
	\end{eqnarray*}
Recalling the Riemannian metric \eqref{eq:riemannmetricF}, points to the steepest ascent direction
	 \[
	 	v=-\nabla_L \log\rho,~~ s=-\frac{1}{\alpha} \log\rho.
	 \]
Thus, the gradient flow of $\cS$ with respect to $W_{2,F}$ is
	\[
		\dot\rho(t)=-\frac{1}{2} \nabla_L^* (\rho \nabla_L\log\rho+\nabla_L\log\rho\cdot\rho)
		-\frac{1}{\alpha}\log\rho.
	\]	
Similarly, for $\cU$ we have
	\begin{eqnarray*}
	\frac{d\cU(\rho(t))}{dt}&=& \trace((\rho-\hat\rho)\dot\rho)
	\\&=&\trace((\rho-\hat\rho)(\frac{1}{2} \nabla_L^* (\rho v+v\rho)+s))
	\\&=&\trace(\frac12\rho((\nabla_L(\rho-\hat\rho))^*v+v^*\nabla_L(\rho-\hat\rho))+(\rho-\hat\rho) s),
	\end{eqnarray*}
which gives steepest descent direction
	\[
		v=-\nabla_L(\rho-\hat\rho),~~s=-\frac{1}{\alpha}(\rho-\hat\rho),
	\]
and the gradient flow
	\[
		\dot\rho(t)=-\frac{1}{2} \nabla_L^* (\rho \nabla_L(\rho-\hat\rho)+\nabla_L(\rho-\hat\rho)\rho)
		-\frac{1}{\alpha}(\rho-\hat\rho).
	\]

\section{Conclusions}

Our line of research into unbalanced versions of optimal mass transport is motivated by the fact that general distributions (histograms, power spectra, spatio-temporal energy densities, images) may not necessarily be normalized to have the same integral. Thus, it is imperative to devise appropriate metrics and theory to handle these situations.
Our overall aim is to provide constructions for ``interpolating'' data in the form of distributions. In the present work, we have formulated a natural technique that interpolates the quantum mechanical version of OMT developed in \cite{Chen} with an analogue of Fisher-Rao information.

In further work, we plan to explore the associated Riemannian structure associated to the unbalanced Wasserstein distance, gradient flow of entropy, and other variants of the continuity equation. From a more applied side, we plan to implement the methodology described in the present work in Matlab (as noted above it can be numerically solved via convex optimization), and then test it various types of multi-modal, multi-sensor, and multi-spectral data. It seems ideal for multiple target estimation as was done in \cite{Lipeng}.

\section*{Acknowledgements}
This project was supported by AFOSR grants (FA9550-15-1-0045 and FA9550-17-1-0435), grants from the National Center for Research Resources (P41-
RR-013218) and the National Institute of Biomedical Imaging and Bioengineering (P41-EB-015902), National Science Foundation (NSF), and a postdoctoral fellowship through Memorial Sloan Kettering Cancer Center.

\bibliographystyle{plain}

\end{document}